\documentclass[12pt,eqright]{amsart}
\usepackage{amssymb,amsmath}
\usepackage{amsfonts}
\usepackage{epsf}
\usepackage{graphicx}
\newtheorem{theorem}{Theorem}[section]
\newtheorem{lemma}[theorem]{Lemma}

\newtheorem{remark}[theorem]{Remark}

\theoremstyle{definition}
\newtheorem{algorithm}{Algorithm}[section]

\numberwithin{equation}{section}

\newcommand{\bc}{\begin{center}}
\newcommand{\ec}{\end{center}}
\newcommand{\be}{\begin{eqnarray}}
\newcommand{\ee}{\end{eqnarray}}

\newcommand{\ben}{\begin{eqnarray*}}
\newcommand{\een}{\end{eqnarray*}}

\newcommand{\Om}{\Omega}
\newcommand{\al}{\alpha}

\newcommand{\pa}{\partial}

\newcommand{\na}{\nabla}

\def\x{\times}

\def\na{\nabla}

\def\ds{\,ds}

\def\cC{\mathcal{C}}

\def\cE{\mathcal{E}}

\def\cT{\mathcal{T}}
\def\R{\mathbb{R}}

 \DeclareMathOperator{\res}{Res}
 \DeclareMathOperator{\tr}{tr}

\DeclareMathOperator{\card}{card}
\DeclareMathOperator{\osc}{osc}

\DeclareMathOperator{\consis}{\kappa}

\title[Convergence and optimality of adaptive nonconforming methods]
{Convergence and optimality of the adaptive Morley element method}
\author[J.~Hu]{Jun Hu}
\address{LMAM and  School of Mathematical Sciences,
  Peking University,  Beijing 100871, P. R. China}
\email{hujun@math.pku.edu.cn}

\author[Z.~C.~Shi]{Zhongci Shi}
\address{Institute of Computational Mathematics, Chinese Academy of Sciences, Beijing 100080, P.~R.~China.}
\email{shi@lsec.cc.ac.cn}

\author[J.~Xu]{Jinchao Xu}
\address{The School of Mathematical Sciences, Peking University, and
  Department of Mathematics, Pennsylvania State University, University
  Park, PA 16801}
\email{xu@math.psu.edu}

\thanks{The first author was supported by the  NSFC under Grant 10971005, and
 Foundation for the Author of National Excellent Doctoral
Dissertation of PR China 200718, and partially supported by  the Chinesisch-Deutsches Zentrum project
GZ578. The
third author was supported in part by NSF 0915153  and
NSFC-10528102.}

\begin{document}
\begin{abstract}
  This paper is devoted to the convergence and optimality analysis of
   the adaptive Morley element method for  the fourth order elliptic problem.
  A new technique is developed to establish  a quasi-orthogonality
  which is crucial for the convergence analysis of the adaptive
  nonconforming method.  By introducing a new  parameter-dependent error estimator and further establishing a
  discrete reliability property, sharp convergence and optimality
  estimates are then fully proved for the fourth order elliptic problem.  
   
\end{abstract}
\maketitle
\section{Introduction}
  This paper is devoted to the study of
adaptive nonconforming  finite element methods for high order
elliptic boundary value problems. The adaptive conforming finite
element method for the second order elliptic problems has been a
subject of extensive studies for many years since the pioneering
work of Babuska and Rheinboldt \cite{BaRe78}, and its theory  has
become rather mature
\cite{Verfurth96,Dolfler96,MorinNochettoSiebert00,MorinNochettoSiebert02,
  BinevDahmenDeVore04,Stevenson05,Stevenson06,CasconKreuzerNochettoSiebert07}.
For the nonconforming method, the a posteriori error theory of the
second order elliptic problems has been studied only very recently
\cite{DarDurPad95,DarDurPadVam96,CarBarJan02,C,CarstensenHuOrlando07,CarstensenHu07};
for the fourth order elliptic problem,  only the a posteriori error
estimate of the Morley element method can be found in the literature
\cite{HuShi08,St07,ZW06} and there have been no works on either
convergence or optimality for any finite element methods for
fourth order problems.

 The main difficulty for the analysis of nonconforming finite element methods
  arises  from the nonconformity of
the discrete space and consequently the lack of the
Galerkin-orthogonality which is a key ingredient for the convergence
analysis of the adaptive conforming method of the second order elliptic
problem
\cite{Dolfler96,MorinNochettoSiebert00,MorinNochettoSiebert02,Stevenson06,CasconKreuzerNochettoSiebert07}.
For the nonconforming linear element of  the Poisson equation, a
quasi-orthogonality is established instead in \cite{CarstensenHoppe05b} by
using  some special equivalency between the nonconforming linear
element and the lowest order Raviart-Thomas element \cite{Marini85}.
For the Morley element of the fourth order elliptic problem,
however, it is  unclear whether  such type of
equivalency  still holds.  We also note that  the convergence (not
to mention optimality) analysis of the adaptive conforming method is
still missing  for the fourth order elliptic problem in the
literature.

This paper is devoted to the convergence and optimality  analysis of
the adaptive version of the Morley element \cite{Morley68,
Shi90,WX06}. Our analysis  is based on an observation that a
quasi-orthogonality can be obtained from a crucial local
conservative property (that plays a critical role in a general study
in \cite{WX06}), of the  Morley  element method. Another ingredient
is  a new parameter dependent estimator which is introduced to
analyze  optimality of the adaptive nonconforming method. With
the help of the discrete reliability which is established by
introducing two interpolation operators between two nonconforming
spaces, we show  convergence and optimality of the adaptive
algorithm.

 The rest of the paper is organized as follows. In Section 2, we
 present the Kirchhoff plate problem and the Morley finite element
 method, and recall a posteriori error analysis  due to \cite{HuShi08}. In
 Section 3, we prove the quasi-orthogonality and then  show  reduction
 of some total error in Section 4  by introducing a new
 parameter-dependent estimator. To obtain  optimality of the
 adaptive algorithm,  we establish the discrete reliability in
 Section 5. Consequently, we show  optimality of the adaptive Morley
 element method in Section 6.  We give a brief comment on the extension of the theory
 to the Morley element method in three dimensions in Section 7.  Also, we
 discuss the generalization to the nonconforming linear
 elements in both two and three dimensions therein. This
 extension  gives an alternative  analysis of the convergence result from
 \cite{CarstensenHoppe05b}.  The paper ends with Section 8 where we give the conclusion and some comments.

\section{The Morley element for the Kirchhoff plate problem and an a posteriori error estimate}
Let $\Om\subset \R^2$ be a bounded  domain, $\mathbb{E}$  the Young
modulus, and $\nu$ the Poisson ratio.  For all $2\x 2$ symmetric
matrices, the linear operator $\cC$ is defined by
$$
\cC\tau:=\frac{\mathbb{E}}{12(1-\nu^2)}\big((1-\nu)\tau+\nu\tr(\tau)I\big)\,.
$$
The  bilinear form $a(u, v)$  is defined by
 \begin{equation}\label{vararm}
  a(u,v)=(\cC\na^2 u,\na^2 v)_{L^2(\Om)}, \text{ for any } u\,,
 v\in W:= H_0^2(\Om)\,,
\end{equation}
where $\na^2 u$ is the Hessian matrix of $u$.  The corresponding  energy norm is given by
\begin{equation}
\|u\|_{\cC}^2:=a(u, u) \text{ for any }
 u\in W\,,
\end{equation}
which is equivalent to the usual norm $|\cdot|_{H^2(\Om)}$ for any
$u\in W$.

We consider the Kirchhoff plate bending problem as follows: Given
$f\in L^2(\Om)$,  find $u\in W$ such that
\begin{equation} \label{rmplate}
 a(u, v)=(f,v)_{L^2(\Om)} \quad \text{ for all }  v\in W\,.
\end{equation}

We now present the Morley element.  Suppose
that  $\overline{\Omega}$ is covered exactly by a shape-regular
triangulation $\cT_h$ consisting of
 triangles in $2D$, see \cite{CiaBook}.  $\cE_h$ is the set of all  edges in $\cT_h$,
$\cE_h(\Om)$ is the set of interior edges, and $\cE(K)$ is the set of
edges of any given element $K$ in $\cT_h$; $h_K=|K|^{1/2}$,
 the size of the element $K\in \cT_h$.
$\omega_K$ is the union of elements $K'\in \cT_h$  that  share an
edge with $K$, and $\omega_e$ is  the union of elements that share a
common edge $e$. Given any edge $e\in\cE_h(\Om)$ with the length $h_e$
we assign one fixed unit normal $\nu_e:=(\nu_1,\, \nu_2)$ and
tangential vector $\tau_e:=(-\nu_2,\,\nu_1)$. For $e$ on the
boundary we choose $\nu_e=\nu$ the unit outward normal to $\Omega$.
Once $\nu_e$ and $\tau_e$ have been fixed on $e$, in relation to
$\nu_e$ one defines the  elements $K_{-}\in \cT_h$ and $K_{+}\in
\cT_h$, with $e=K_{+}\cap K_{-}$.  Given $e\in\cE_h(\Omega)$ and some
$\R^d$-valued function $v$ defined in $\Omega$, with $d=1,2$, we
denote by $[v]:=(v|_{K_+})|_e-(v|_{K_-})|_e$ the jump of $v$ across
$e$.

The discrete space of the  Morley finite element method  is defined
as follows \cite{Morley68,Shi90,WX06}
\begin{equation}
\begin{split}
W_h := &\{v\in M_{2,h}, \int_e[\na_h v \cdot \nu_e]ds=0 \text{
on } e\in\cE_h(\Om),\\
&\text{ and } \int_e \na v\cdot \nu_eds=0 \text{ on }
e\in\cE_h\cap\pa\Om\} \mbox{,}
\end{split}
\end{equation}
where $M_{2,h}$ is the space of piecewise polynomials of degree
$\leq$ 2 over $\cT_h$ which are continuous at all the internal nodes
and vanish at all the nodes on the boundary $\pa\Om$, and  $\na_h$
the discrete gradient operator which is defined elementwise. We
define
\begin{equation}\label{discreteenergy}
\begin{split}
a_h(u_h, v_h):&=(\cC\na_h^2 u_h,\na_h^2 v_h)_{L^2(\Om)}
\text{ for any }  u_h\,, v_h\in  W+W_h\,,\\[0.5ex]
\|u_h\|_{\cC_h}^2:&=a_h(u_h, u_h) \text{ for any }
 u_h\in W+W_h\,,
\end{split}
\end{equation}
where the discrete Hessian operator $\na_h^2$ is defined elementwise with respect to the
triangulation $\cT_h$.

We now consider the finite element discretization of \eqref{rmplate}
as follows: Find $ u_h\in W_h$ such that
 \begin{equation}\label{discrm}
  a_h(u_h,v_h) =(f,v_h)_{L^2(\Om)} \quad\text{ for all } v_h\in W_{h}\,.
 \end{equation}

To  recall  the a posteriori error estimate for the Morley element, we first
define an estimator on each element $K\in\cT_h$ as
\begin{equation}\label{estimator}
\eta_K=h_K^{2}\|f\|_{L^2(K)}+\bigg(\sum\limits_{ e\subset\pa
K}h_K\|[\na_h^2 u_h\tau_e]\|_{L^2(e)}^2\bigg)^{1/2}\,.
\end{equation}
For any $S_h\subset\cT_h$, we define the estimator over $S_h$  by
\begin{equation}\label{eq2.6a}
\eta^2(u_h, S_h):=\sum\limits_{K\in S_h}\eta_K^2\,.
\end{equation}
In particular, for $S_h=\cT_h$, we have
\begin{equation}\label{eq2.6}
\eta^2(u_h, \cT_h):=\sum\limits_{K\in \cT_h}\eta_K^2\,.
\end{equation}
 We further define the
 oscillation $\osc(f, \cT_h)$ by
\begin{equation}\label{eq2.8}
\osc^2(f,\cT_h):=\sum\limits_{K\in\cT_h}h_K^{4}\|f-f_K\|_{L^2(K)}^2\,,
\end{equation}
where $f_K$ is  the constant projection of $f$ over $K$. For the
estimator \eqref{eq2.6}, we have the following reliability and
efficiency whose proof can be found in \cite{HuShi08}.
\begin{lemma}\label{Lemma2.1} Let $u$  be the solution of Problem \eqref{rmplate},
and $u_h$ be the solution of  Problem \eqref{discrm}. Then,
\begin{equation}\label{eq2.7}
\begin{split}
\|u-u_h\|_{\cC_h} &\lesssim \eta(u_h, \cT_h)
\lesssim \|u-u_h\|_{\cC_h}+\osc(f,\cT_h)\,.
\end{split}
\end{equation}
\end{lemma}
Here and throughout the paper, we shall follow \cite{Xu92} to use
the notation $\lesssim$ and $\approxeq$. When we write
$$
A_1 \lesssim B_1, \text{ and } A_2\approxeq B_2,
$$
then there exist possible constants $C_1$, $c_2$ and $C_2$ such that
$$
A_1 \leq C_1 B_1, \text{ and } c_2B_2\leq A_2\leq C_2 B_2.
$$

Given $v\in H^2(\cT_h):=\{v\in L^2(\Om), v|_K\in H^2(K), \text{ for
any }K\in\cT_h\}$,  we define the following residual
\begin{equation}\label{eq2.9}
\res_H(v)=(f, v)_{L^2(\Om)}-a_h(u_H,v), \text{ for any }v\in
H^2(\cT_h)\,,
\end{equation}
with $u_H$ being the solution of the discrete problem \eqref{discrm}
on  $\cT_H$, which is a nested and coarser mesh to $\cT_h$;
namely, $\cT_h$ is some refinement of $\cT_H$. It follows from the
discrete problem \eqref{discrm} that
\begin{equation}\label{eq2.10}
\res_H(v)=\res_H(v-v_H), \text{ for any }v_H\in W_H.
\end{equation}

\begin{lemma}  For any $v\in W$, it  holds that
\begin{equation}
|\res_H(v)|\lesssim
\big(\sum\limits_{K\in\cT_H}h_K^{4}\|f\|_{L^2(K)}^2\big)^{1/2}\|v\|_{\cC}
\text{ for any }
 v\in W.
\end{equation}
\end{lemma}
 The proof of the above lemma can be found  in  \cite{ZW06,St07,HuShi08}.
\qed

\section{Quasi-orthogonality}
In this section, we address one difficulty, namely
 the quasi-orthogonality,  in the convergence analysis of the
adaptive Morley element method. Our analysis is based on two interpolation
operators: the canonical interpolation operator $\Pi_h$ of the
nonconforming space $W_h$,  and the restriction operator $I_H$ from
the discrete space $W_h$ on the mesh $\cT_h$ to the discrete space
$W_H$ on the nested and coarser mesh $\cT_H$ of $\cT_h$.

Here and in what follows, $\mathcal{N}_h$ denotes the set of nodes
of the partition $\cT_h$.  We first define the canonical
interpolation operator $\Pi_h: W\rightarrow W_h$ by,
\begin{equation}\label{eq3.1}
\begin{split}
(\Pi_hv)(P)=v(P),  \int_e\na_h (\Pi_hv-v)\cdot \nu_e\ds=0 \text{ for
any } v\in W,  P\in\mathcal{N}_h, e\in\cE_h\,.
\end{split}
\end{equation}
\begin{lemma} Let the interpolation operator $\Pi_{h}$  be defined
as in \eqref{eq3.1}. Then,
\begin{equation}\label{eq3.2a}
\int_e \nabla_h (v-\Pi_hv)\ds=0 \text{ for any }e\in\cE_h\text{ and
}v\in W\,,
\end{equation}
\begin{equation}\label{eq3.2ab}
a_h(v-\Pi_hv,v_h)=0\text{ for any } v\in W, v_h\text{ piecewise
quadratic }\,,
\end{equation}
\begin{equation}\label{eq3.4a}
\|v-\Pi_hv\|_{L^2(K)}\lesssim h_K^2|v|_{H^2(K)}\text{ for any
}K\in\cT_h\text{ and }v\in W\,.
\end{equation}
\end{lemma}
The above properties are immediate from the definition of $\Pi_h$.
Now we define
 the restriction interpolation operator $I_H: W_h\rightarrow W_H$ by, for any $v_h\in W_h$,
\begin{equation}\label{eq3.3}
\left\{
\begin{split}
 &(I_Hv_h)(P)=v_h(P), P\in\mathcal{N}_H\,,\\[0.5ex]
 &\int_e \frac{\pa (I_Hv_h)}{\pa
\nu_{e}}ds=\sum\limits_{l=1}^{\ell}\int_{e_l} \frac{\pa v_h}{\pa
\nu_{e} }ds\,,  e\in \cE_{H} \text{ with }e=e_1\cup e_2\cdots\cup
e_{\ell} \text{ and }e_i\in\cE_{h} \,.
\end{split}
\right.
\end{equation}

Before analyzing the properties of this interpolation, we  state the
following simple result.
\begin{lemma}\label{averagecon}  Let $K_1, K_2\in\cT_h$ be two elements sharing a common edge $e$.
If $v_h\in W_h(K_1\cup K_2)$ and $\na^2_hv_h=0$, then $v_h\in
P_1(K_1\cup K_2)$. Namely $v_h$ is a polynomial of degree $\leq 1$
over $K_1\cup K_2$.
\end{lemma}
\begin{proof}
By the definition of $W_h$, $v_h$ is continuous on $K_1\cup K_2$.
Further $\frac{\pa v_h}{\pa\nu_e}|_{K_1}$ and $\frac{\pa
v_h}{\pa\nu_e}|_{K_2}$ are two constant functions that must be equal
since by the definition of $W_h$ $\int_e[\frac{\pa
v_h}{\pa\nu_e}]ds=0$. Thus $v$ must belong to $P_1(K_1\cup K_2)$.
\end{proof}

The properties of the interpolation operator $I_H$ are summarized in
the following lemma.
\begin{lemma} Let the interpolation operator $I_H$  be defined
as in  \eqref{eq3.3}. Then,
\begin{equation}\label{eq3.2b}
\int_e \nabla_h (v_h-I_Hv_h)\ds=0 \text{ for any }e\in\cE_{H}\text{ and
}v_h\in W_h\,,
\end{equation}
\begin{equation}\label{eq3.2bb}
a_h(v_H, v_h-I_Hv_h)=0\text{ for any }v_H\in W_H, v_h\in W_h\,,
\end{equation}
\begin{equation}\label{eq3.4b}
I_Hv_h|_K=v_h|_K  \text{ for any } K\in \cT_h\cap\cT_H\text{ and }v_h\in W_h\,,
\end{equation}
\begin{equation}\label{eq3.5b}
\|I_Hv_h-v_h\|_{L^2(K)}\lesssim h_K^2\|\nabla_h^2v_h\|_{L^2(K)} \text{ for
any }K\in\cT_H\backslash\cT_h\text{ and }v_h\in W_h\,.
\end{equation}
\end{lemma}
\begin{proof} The properties of \eqref{eq3.2b}, and
\eqref{eq3.4b} directly follow  from the definition of the
interpolation. We only need to prove  \eqref{eq3.2bb} and the
estimate \eqref{eq3.5b}.

We first define $\sigma_H=\cC\na_H^2 v_H$ to  assert that
\begin{equation}\label{eq3.10}
\int_e \nabla_h (I-I_H)v_h\cdot \sigma_H\nu_e\ds=0 \text{ for any
}e\in\cE_{H}\,.
\end{equation}
In fact, for $e\in\cE_{H}\backslash\cE_{h}$, this assertion follows
from the fact that $\sigma_H$ is a piecewise constant matrix  with respect to
$\cT_H$ and the definition of $I_H$ in \eqref{eq3.3}.
 For
$e\in\cE_{h}\cap\cE_{H}$, the assertion follows from
$(I-I_H)v_h|_e=0$.

 For the edge $e\in\cE_{h}$ which lies in the
interior of some  $K\in\cT_H$,  we can use the continuity of $\int_e
\nabla_hv_h\ds$ over $e$ and the fact $\sigma_H$ is constant over
$K$ to show that
\begin{equation}\label{eq3.10b}
\int_e[\nabla_h (I-I_H)v_h]\cdot \sigma_H\nu_e\ds=0\,.
\end{equation}
Whence, we integrate by parts and use \eqref{eq3.10} and
\eqref{eq3.10b} to conclude \eqref{eq3.2bb}.

Now we turn to \eqref{eq3.5b}.  In fact, both sides of
\eqref{eq3.5b} are semi-norms of the restriction $W_h(K)$ of $W_h$
on $K$. If the right hand side vanishes for some $v_h\in W_h(K)$, then
$v_h$ is a piecewise polynomial of degree$\leq 1$ on $K$ with respect
to $\cT_h$. It follows from Lemma \ref{averagecon} that $v_h$ is a
polynomial of degree$\leq 1$ on $K$. Therefore the left hand side
also vanishes for the same $v_h$. The desired result then follows from
a scaling argument.
\end{proof}

\begin{lemma}\label{Lemma3.2}(Quasi-orthogonality) Let $\cT_h$ be a refinement of
$\cT_H$,  and $u_h$ and $u_H$ be the solutions of \eqref{discrm} on
$\cT_h$ and $\cT_H$, respectively. Then,
\begin{equation}\label{eq3.6}
\begin{split}
&|a_h(u_h-u_H, u-u_h)|\lesssim
\sum\limits_{K\in\cT_H\backslash\cT_h}h_K^2\|f\|_{L^2(K)}\|\nabla_h^2(u-u_h)\|_{L^2(K)}\,.
\end{split}
\end{equation}
\end{lemma}
\begin{proof} Let the interpolation operator $\Pi_h$ be defined as in
\eqref{eq3.1}.  Since $\Pi_h$ is well-defined for any $v_h\in W_h$(
in fact, $\Pi_hv_h=v_h$) and $a_h(u_h-u_H, (I-\Pi_h)(u-u_h))=0$ (by
\eqref{eq3.2ab}), we have
\begin{equation}\label{eq3.7}
a_h(u_h-u_H, u-u_h)=a_h(u_h-u_H,\Pi_h(u-u_h)).
\end{equation}
Let $v_h=\Pi_h(u-u_h)$ and  the interpolation $I_Hv_h$ be defined as
in \eqref{eq3.3}. The combination of  \eqref{discrm} and \eqref{eq2.9} leads to
\begin{equation}\label{eq3.9}
\begin{split}
a_h(u_h-u_H, v_h)=(f,
v_h)_{L^2(\Om)}-a_h(u_H,v_h)\\
=(f, (I-I_H)v_h)_{L^2(\Om)}-a_h(u_H,(I-I_H)v_h)\,.
\end{split}
\end{equation}
By \eqref{eq3.4b} and \eqref{eq3.5b}, we have
\begin{equation}
|(f, (I-I_H)v_h)_{L^2(\Om)}| \lesssim
\sum\limits_{K\in\cT_H\backslash\cT_h}h_K^2\|f\|_{L^2(K)}\|\nabla_h^2
v_h\|_{L^2(K)}\,.
\end{equation}
From \eqref{eq3.2bb} we have $a_h(u_H, (I-I_H)v_h)=0$.  Then, the
desired result then follows from the triangle inequality and the
approximation property of $\Pi_h$.
\end{proof}

\begin{remark}
For the nonconforming $P_1$ element of the Poisson equation, the
quasi-orthogonality  was  established in \cite{CarstensenHoppe05b}.
The analysis therein is based on some special equivalency between
the nonconforming $P_1$ and Raviart-Thomas elements.  For the
Stokes-like problem, the quasi-orthogonality of the nonconforming
$P_1$ element has been first proved in \cite{HuXu2007} based on some
special relation of the nonconforming $P_1$ and Raviart-Thomas
elements. For the Morley element, it is unclear whether there exists
similar equivalency or relation so far.
\end{remark}

\begin{remark}
 This paper is a refined version of  a  technical report in 2009 \cite{HuShiXu2009}, where it was the first time in the literature to
  make use of the conservative property of the nonconforming finite element space to analyze the quasi-orthogonality.
\end{remark}

\section{Reduction of a properly defined total error}
In the rest of the paper, we shall establish  convergence and
optimality of our Adaptive Nonconforming Finite Element
Method(ANFEM). Our analysis is based on two main ingredients: the
strict reduction of some total error between two levels and the
discrete reliability of the estimator.  To this end, we shall first
introduce a modified estimator $ \tilde{\eta}$ with a
 undetermined  positive constant; we shall then borrow the concept
of the total error of
\cite{CasconKreuzerNochettoSiebert07,CarstensenHoppe05b} which
contains the energy norm of the error and the scaled  estimator
$\tilde{\eta}$;  we  shall finally show  reduction of this total
error.  We shall establish  the discrete reliability of the
estimator in the next section.

Let us first define our adaptive algorithm.  Given an initial shape
regular triangulation $\mathcal T_0$, a right-hand side function
$f\in L^2(\Omega)$, a tolerance $\varepsilon>0$, and a parameter
$\theta\in(0,1)$. Hereafter, we shall replace the subscript $h$ by
an iteration counter called $k$.
\begin{algorithm}\label{Algorithm}
\noindent $[\mathcal{T}_N\,, u_N]$={\bf \small
ANFEM}$(\mathcal{T}_0, f,\varepsilon, \theta)$

$\eta = \varepsilon\,, k=0$

\smallskip

{\bf \small WHILE} $\eta \geq \varepsilon$, {\bf \small DO}

(1)\quad Solve \eqref{discrm} on $\mathcal{T}_k$,  to get the solution
$u_k$.

(2)\quad Compute the error estimator $\eta=\eta(u_k, \cT_k)$.

(3)\quad Mark the minimal element set $\mathcal{M}_k$ such that
  \begin{equation}\label{eq4.1}
 \eta^2(u_k, \mathcal{M}_k)\geq \theta \, \eta^2 (u_k, \cT_k).
  \end{equation}

(4)\quad Refine each triangle $K \in \mathcal{M}_k$
  by the newest vertex bisection and   possible further refining to conformity to get $\cT_{k+1}$.

\quad $k=k+1$.

\smallskip

{\bf \small END WHILE}

\smallskip

$\cT_N=\cT_k$.

\noindent {\bf \small END ANFEM}
\end{algorithm}

In order to  prove  a strict reduction of some total error, we
 define the following modified estimator
\begin{equation}\label{eq5.2}
\tilde{\eta}^2(u_H, \cT_H):=
\sum\limits_{K\in\cT_H}\big(\beta_1h_K^{4}\|f\|_{L^2(K)}^2
 +\eta_K^2\big) \text{ with
}\eta_K \text{ defined in } \eqref{estimator}
\end{equation}
for some positive constant $\beta_1$ to be determined later.
\begin{remark}Note that, as we can see below,  the modified error estimator $\tilde{\eta}(u_H,
\cT_H)$ is only for the analysis,  the final results concerning
both convergence and   optimality will be  proved for Algorithm
\ref{Algorithm}.
\end{remark}

\begin{lemma}\label{Lemma5.3b} Let $\cT_h$ be some refinement of $\cT_H$ with the bulk  criterion \eqref{eq4.1},
 then there exist $\rho>0$ and a
  positive constant $\beta\in(1-\rho\theta,1)$
such that
\begin{equation}\label{eq5.7}
{\eta}^2(u_H,\cT_{h})
 \leq  \beta {\eta}^2(u_H,\cT_{H})+(1-\rho\theta-\beta){\eta}^2(u_H,\cT_{H})\,.
\end{equation}\end{lemma}
\begin{proof} The result can be proved by following the idea in \cite{CasconKreuzerNochettoSiebert07}.
 We give the details only for the readers' convenience.  In fact, we have
\begin{equation}
{\eta}^2(u_H,\cT_{h})={\eta}^2(u_H,\cT_{H}\cap\cT_h)+{\eta}^2(u_H,\cT_{H}\backslash\cT_h).
\end{equation}
For any $K\in \cT_H\backslash\cT_h$, we only need to consider the
case where  $K$ is subdivided into $K_{1}\,,K_2\in\cT_h$ with
$|K_1|=|K_2|=\frac{1}{2}|K|$. By the definitions of $h_K$ and
${\eta}_K(u_H)$, we have
\begin{equation}\label{eq4.5}
\begin{split}
 &\sum\limits_{i=1}^2{\eta}_{K_i}^2(u_H)
 :=\sum\limits_{i=1}^2\bigg(h_{K_i}^{2}\|f\|_{L^2(K_i)}+\bigg(\sum\limits_{\cE_{h}\ni
e\subset\pa
K_i}h_{K_i}\|[\nabla_H^2u_H \tau_e]\|_{L^2(e)}^2\bigg)^{1/2}\bigg)^2\\[0.5ex]
& \leq\frac{1}{2^{1/2}} {\eta}_K^2(u_H):=\frac{1}{2^{1/2}}
 \bigg( h_{K}^{2}\|f\|_{L^2(K)}+\bigg(\sum\limits_{\cE_{H}\ni
e\subset\pa
K}h_{K}\|[\nabla_H^2u_H\tau_e]\|_{L^2(e)}^2\bigg)^{1/2}\bigg)^2\,,
 \end{split}
\end{equation}
since $[\nabla_H^2u_H \tau_e]=0$ over $e=K_1\cap K_2\in\cE_h$.
Consequently
\begin{equation}
\begin{split}
\sum\limits_{K\in\cT_H\backslash\cT_h}
\sum\limits_{i=1}^2{\eta}_{K_i}^2(u_H) \leq
\frac{1}{2^{1/2}}{\eta}^2(u_H, \cT_H\backslash\cT_h)\,,
\end{split}
\end{equation}
and
\begin{equation}
\begin{split}
{\eta}^2(u_H,\cT_{h}) &\leq {\eta}^2(u_H,\cT_{H})-\rho{\eta}^2(u_H,
\cT_H\backslash\cT_{h})\,,
\end{split}
\end{equation}
with $\rho=1-\frac{1}{2^{1/2}}$. Taking the positive parameter
$\beta$ with $1-\rho\theta<\beta<1$, the desired result follows by
combining the above inequality and the bulk criterion \eqref{eq4.1}.
\end{proof}

\begin{lemma}
Let $\cT_h$ be some refinement of $\cT_H$,
 then there exists $\rho>0$  such that
 \begin{equation}\label{righthandside}
 \sum\limits_{K\in\cT_h}h_K^{4}\|f\|_{L^2(K)}^2\leq
 \sum\limits_{K\in\cT_H}h_K^{4}\|f\|_{L^2(K)}^2-\rho\sum\limits_{K\in\cT_H\backslash\cT_h}h_K^{4}\|f\|_{L^2(K)}^2\,.
 \end{equation}
\end{lemma}
\begin{proof}  The proof immediately follows from the definition of
 the meshsize $h_K$.
\end{proof}

\begin{lemma}(Continuity of the estimator)  Let $u_h$ and $u_H$ be the solutions to the discrete problem \eqref{discrm}
on the meshes $\cT_h$ and $\cT_H$, respectively. Given any positive
constant $\epsilon$, there exists a positive constant
$\beta_2(\epsilon)$ dependent on $\epsilon$ such that
\begin{equation}\label{eq5.4}
{\eta}^2(u_{h}, \cT_{h})\leq (1+\epsilon) {\eta}^2(u_H, \cT_{h})
+\frac{1}{\beta_2(\epsilon)}\|u_{h}-u_H\|_{\cC_{h}}^2\,.
\end{equation}
\end{lemma}
\begin{proof}
Given any $K\in\cT_h$,  it follows from the definitions of
${\eta}_K(u_h)$ and ${\eta}_K(u_H)$ in \eqref{eq4.5} that
\begin{equation}
\begin{split}
\big |{\eta}_K(u_h)-{\eta}_K(u_H)\big |&=\bigg
|\bigg(\sum\limits_{\cE_{h}\ni e \subset\pa K}h_K\|[\nabla_h^2
 u_h\tau_e]\|_{L^2(e)}^2\bigg)^{1/2}\\
 &\quad -\bigg(\sum\limits_{\cE_{h}\ni
e \subset\pa K}h_K\|[\nabla_H^2
 u_H\tau_e]\|_{L^2(e)}^2\bigg)^{1/2}\bigg |\\
 &\leq \bigg(\sum\limits_{\cE_{h}\ni
e \subset\pa K}h_K\|[\nabla_h^2
 (u_h-u_H)\tau_e]\|_{L^2(e)}^2 \bigg)^{1/2}.
 \end{split}
\end{equation}
With $e=K_1\cap K_2\in\cE_h$,  we use the trace theorem and the fact
 that $\nabla_h^2(u_h-u_H)$ is a piecewise constant matrix to get
\begin{equation}
\begin{split}
\|[\nabla_h^2
 (u_h-u_H)\tau_e]\|_{L^2(e)}
 &\leq \|\nabla_h^2
 (u_h-u_H)\tau_e|_{K_1}\|_{L^2(e)}+\|\nabla_h^2
 (u_h-u_H)\tau_e|_{K_2}\|_{L^2(e)}\\[0.5ex]
 &\lesssim h_K^{-1/2}\|\nabla_h^2
 (u_h-u_H)\|_{L^2(\omega_e)}\,,
 \end{split}
\end{equation}
which gives
\begin{equation}
\big |{\eta}_K(u_h)-{\eta}_K(u_H)\big | \lesssim \|\nabla_h^2
 (u_h-u_H)\|_{L^2(\omega_K)}.
\end{equation}
Applying the Young inequality with any positive constant $\epsilon$
and summarizing over all elements in $\cT_h$ completes the proof of
the lemma.
\end{proof}

\begin{theorem}\label{Theorem5.2} Let $u$ be the solution to the problem \eqref{rmplate},
and $u_H$ and $u_{h}$ be the solutions to the discrete problem
\eqref{discrm} on the meshes $\cT_H$ and $\cT_{h}$, respectively.
Then, there exists positive constants $\gamma_1$, $\beta_1$, and
$0<\al<1$ with
\begin{equation}\label{eq5.3}
\|u-u_{h}\|_{\cC_{h}}^2+\gamma_1\tilde{\eta}^2(u_{h}, \cT_{h})\leq
\al(\|u-u_{H}\|_{\cC_{H}}^2+\gamma_1\tilde{\eta}^2(u_{H},
\cT_{H}))\,.
\end{equation}
\end{theorem}
\begin{proof}   Let
$\delta$, $\gamma_1$,  and $\gamma_2$,   be three positive constants
to be chosen later. Applying the Young inequality to Lemma
\ref{Lemma3.2} and adding the resulting estimate to the inequality
\eqref{eq5.4} leads to
\begin{equation}\label{eq5.5}
\begin{split}
&(1-\delta)\|u-u_{h}\|_{\cC_{h}}^2+\gamma_1
{\eta}^2(u_{h}, \cT_{h})+\gamma_2\sum\limits_{K\in\cT_h}h_K^{4}\|f\|_{L^2(K)}^2\\[0.5ex]
&\leq \|u-u_{H}\|_{\cC_H}^2
 +\gamma_1(1+\epsilon){\eta}^2(u_{H},\cT_{h})
  +(\frac{\gamma_1}{\beta_2(\epsilon)}-1)\|u_{h}-u_H\|_{\cC_{h}}^2\\[0.5ex]
 &\quad +\gamma_2\sum\limits_{K\in\cT_H}h_K^{4}\|f\|_{L^2(K)}^2
 +(C_1(\delta)-\rho\gamma_2)\sum\limits_{K\in\cT_H\backslash\cT_{h}}h_K^{4}\|f\|_{L^2(K)}^2\,,
 \end{split}
\end{equation}
with the positive  constant $\rho$ from \eqref{righthandside}.  We note that the bound for
${\eta}^2(u_{H},\cT_{h})$ is given in Lemma \ref{Lemma5.3b}. Hence
\begin{equation}\label{eq5.5b}
\begin{split}
&(1-\delta)\|u-u_{h}\|_{\cC_{h}}^2+\gamma_1
{\eta}^2(u_{h}, \cT_{h})+\gamma_2\sum\limits_{K\in\cT_h}h_K^{4}\|f\|_{L^2(K)}^2\\[0.5ex]
&\leq
\|u-u_{H}\|_{\cC_H}^2+\gamma_1((1-\rho\theta-\beta)(1+\epsilon)+\epsilon\beta){\eta}^2(u_H,\cT_{H})+\gamma_1\beta
{\eta}^2(u_H,\cT_{H})\\[0.5ex]
&\quad+(\frac{\gamma_1}{\beta_2(\epsilon)}-1)\|u_{h}-u_H\|_{\cC_{h}}^2
+(C_1(\delta)-\rho\gamma_2)
\sum\limits_{K\in\cT_H\backslash\cT_{h}}h_K^{4}\|f\|_{L^2(K)}^2\\[0.5ex]
&\quad+\gamma_2\sum\limits_{K\in\cT_H}h_K^{4}\|f\|_{L^2(K)}^2\,,
\end{split}
\end{equation}
with $\rho$ and $\beta$ from Lemma \ref{Lemma5.3b}. In what follows
we shall choose the parameters $\al$, $\beta$,
 $\gamma_1$, $\gamma_2$, and $\delta$ to achieve the
reduction of the total error. We first set
\begin{equation}
\gamma_2=\frac{C_1(\delta)}{\rho},  \gamma_1=\beta_2(\epsilon),
\text{ and }\beta=(1-\rho\theta)(1+\epsilon)
\end{equation}
which leads to
\begin{equation}\label{eq5.5c}
\begin{split}
&(1-\delta)\|u-u_{h}\|_{\cC_{h}}^2+\gamma_1
{\eta}^2(u_{h}, \cT_{h})+\gamma_2\sum\limits_{K\in\cT_h}h_K^{4}\|f\|_{L^2(K)}^2\\[0.5ex]
&\leq
\|u-u_{H}\|_{\cC_H}^2+\gamma_1\beta{\eta}^2(u_H,\cT_{H})+\gamma_2\sum\limits_{K\in\cT_H}h_K^{4}\|f\|_{L^2(K)}^2\,.
\end{split}
\end{equation}
We choose $\epsilon$ to be small enough such that $0<\beta<1$.  Let
the positive constant $\alpha$ with $\beta<\alpha<1$ be determined
later, this gives
 \begin{equation}\label{eq5.5d}
\begin{split}
&(1-\delta)\|u-u_{h}\|_{\cC_{h}}^2+\gamma_1
{\eta}^2(u_{h}, \cT_{h})+\gamma_2\sum\limits_{K\in\cT_h}h_K^{4}\|f\|_{L^2(K)}^2\\[0.5ex]
&\leq
\alpha\big((1-\delta)\|u-u_{H}\|_{\cC_H}^2+\gamma_1{\eta}^2(u_H,\cT_{H})
+\gamma_2\sum\limits_{K\in\cT_H}h_K^{4}\|f\|_{L^2(K)}^2\big)\\[0.5ex]
&+(1-\alpha(1-\delta))\|u-u_{H}\|_{\cC_H}^2
+\gamma_1(\beta-\alpha){\eta}^2(u_H,\cT_{H})\\[0.5ex]
&+\gamma_2(1-\alpha)\sum\limits_{K\in\cT_H}h_K^{4}\|f\|_{L^2(K)}^2\,.
\end{split}
\end{equation}
 Recalling the reliability of $\eta(u_H, \cT_H)$ with the reliability coefficient $C_{Rel}$ in Lemma
\ref{Lemma2.1}
\begin{equation}\label{eq5.8}
\|u-u_{H}\|_{\cC_H}^2\leq C_{Rel}{\eta}^2(u_H,\cT_{H})\,,
\end{equation}
and the fact that
\begin{equation}
\sum\limits_{K\in\cT_H}h_K^{4}\|f\|_{L^2(K)}^2 \leq
{\eta}^2(u_H,\cT_{H}).
\end{equation}
Whence we derive as
\begin{equation}
\begin{split}
&(1-\alpha(1-\delta))\|u-u_{H}\|_{\cC_H}^2
+\gamma_1(\beta-\alpha){\eta}^2(u_H,\cT_{H})
+\gamma_2(1-\alpha)\sum\limits_{K\in\cT_H}h_K^{4}\|f\|_{L^2(K)}^2
\\[0.5ex]
&\leq
\big((1-\alpha(1-\delta))C_{Rel}+\gamma_1(\beta-\alpha)+\gamma_2(1-\alpha)\big){\eta}^2(u_H,\cT_{H})\,,
\end{split}
\end{equation}
provided that $0<\delta<1$.  Then, the choice of
$\alpha=\frac{\gamma_1\beta+\gamma_2+C_{Rel}}{\gamma_1+\gamma_2+C_{Rel}(1-\delta)}>\beta$
gives
\begin{equation}\label{eq5.5e}
\begin{split}
&(1-\delta)\|u-u_{h}\|_{\cC_{h}}^2+\gamma_1
{\eta}^2(u_{h}, \cT_{h})+\gamma_2\sum\limits_{K\in\cT_h}h_K^{4}\|f\|_{L^2(K)}^2\\[0.5ex]
&\leq
\alpha\big((1-\delta)\|u-u_{H}\|_{\cC_H}^2+\gamma_1{\eta}^2(u_H,\cT_{H})
+\gamma_2\sum\limits_{K\in\cT_H}h_K^{4}\|f\|_{L^2(K)}^2\big)\,.
\end{split}
\end{equation}
We  choose  such that $0<\delta<\min(
\frac{\gamma_1(1-\beta)}{C_{Rel}},1)$ to assure that $\alpha<1$.
Finally, we  take $\beta_1=\gamma_2/\gamma_1$ and redefine
$\gamma_1=\gamma_1/(1-\delta)$ to end the proof.
\end{proof}

\section{Discrete reliability}
This section is devoted to the discrete reliability of the estimator
$\eta(u_H,\cT_H)$.  The analysis needs the prolongation operator
$I_h^{\prime}:W_H\rightarrow W_h$ defined as follows.   Given
$P\in\mathcal{N}_h$ and $e\in\cE_h$, the nodal patch $\omega_{P,H}$
of $P$ and  the edge patch $\omega_{e,H}$ of $e$ with respect to the
mesh $\cT_H$ are defined by, respectively,
\begin{equation}
\begin{split}
 \omega_{P,H}:=\{K\in\cT_H, P\in\partial K \text{ or } P \text{ is in the interior of } K\},\\
  \omega_{e,H}:=\{K\in\cT_H, e\subset\partial K \text{ or }e \text{ is in the interior of }K\}.
 \end{split}
 \end{equation}
Define $\xi_{P}=\card( \omega_{P,H})$ and $\xi_{e}=\card(
\omega_{e,H})$.  We define the prolongation interpolation
$I_{h}^{\prime}v_H\in W_h$ by, for any $v_H\in W_H$,
\begin{equation}\label{eq5.13}
\left\{
\begin{split}
&(I_{h}^{\prime}v_H)(P)=\frac{1}{\xi_P}\sum\limits_{K\in\omega_{P,H}}v_H|_{K}(P)\,\text{ for any }  P\in\mathcal{N}_h\,,\\
& \int_e\frac{\pa (I_{h}^{\prime}v_H)}{\pa\nu_{e}}ds=
\frac{1}{\xi_{e}}\sum\limits_{K\in \omega_{e,H}} \int_e\frac{\pa
(v_H|_{K})}{\pa\nu_{e}}ds\,\text{ for any } e\in\cE_h\,.
\end{split}
\right.
\end{equation}

\begin{lemma}\label{continuity}  Let $K_1, K_2\in\cT_H$ be two elements sharing a common edge $e$
with two endpoints $P_1$ and $P_2$. Suppose that $v_H\in W_H$ and
$\na_Hv_H$ is  continuous over $e$.  Then, $v_H$ is continuous over $e$.
\end{lemma}
\begin{proof}
 We can assume that the common edge $e$ shared by  $K_1$
and $K_2$ lies along the $x$-axis.  Then, $v$ can be expressed as
$$
v_H|_{K_1}=a_0+a_1x+a_2y+a_3xy+a_4x^2+a_5y^2, \text{ and
}v_H|_{K_2}=b_0+b_1x+b_2y+b_3xy+b_4x^2+b_5y^2.
$$
 Since $\frac{\pa v_H}{\pa x}$ is continuous over $e$, we have $a_1=b_1$ and $a_4=b_4$.
 The continuity of $\frac{\pa v_H}{\pa y}$ over $e$ gives $a_2=b_2$ and
 $a_3=b_3$. Finally,  $v_H|_{K_1}(P_{\ell})=v_H|_{K_2}(P_{\ell}),
 \ell=1\,,2$,  concludes $a_0=b_0$.  Therefore, $v_H$ is continuous over $e$.
\end{proof}

\begin{lemma}\label{Lemma5.2} Let the interpolation operator $I_h^{\prime}$  be defined
as in  \eqref{eq5.13}. Then,
\begin{equation}\label{eq3.5cba}
\|\nabla_{h}^2(I_h^{\prime}v_H-v_H)\|_{L^2(\Om)}^2 \lesssim
\sum\limits_{K\in\cT_H\backslash\cT_h}\sum\limits_{e\subset\pa
K}h_K\|[\nabla_H^2 v_H\tau_e]\|_{L^2(e)}^2\,~\text{for any}~ v_H\in W_H.
\end{equation}
\end{lemma}
\begin{proof} It follows from the definition of $I_{h}^{\prime}$
\eqref{eq5.13} that $I_h^{\prime}v_H|_K=v_H|_K$ for any
$K\in\cT_H\cap\cT_h$.  Therefore, we only need to estimate
$\|\nabla_{h}^2(I_h^{\prime}v_H-v_H)\|_{L^2(K)}$ for
$K\in\cT_H\backslash\cT_h$. To prove the desired result, it is
sufficient to show that
\begin{equation}\label{local}
\|\nabla_{h}^2(I_h^{\prime}v_H-v_H)\|_{L^2(K)}\lesssim
\sum\limits_{e\subset\pa K}h_K\|[\nabla_H^2 v_H\tau_e]\|_{L^2(e)}^2
\text{ for any }v_H\in W_H \text{ and }K \in\cT_H\backslash\cT_h.
\end{equation}
For any $e\in\cE_{H}$, $\|[\nabla_H^2v_H\tau_e]\|_{L^2(e)}=0$
indicates that there are no jumps over $e$ for all tangential
components of $\nabla_H^2v_H$, which in turn  implies that $\nabla_Hv_H$
is continuous over $e$ since $\nabla_Hv_H$ is average continuous over
$e$. Since $v_H$ is continuous at all
  the internal nodes,  Lemma \ref{continuity} proves that $v_H$ is
  continuous over $e$. Whence,  $I_h^{\prime}v_H|_K=v_H|_K$
provided that $\|[\nabla_H^2v_H\tau_e]\|_{L^2(e)}=0$ for any $e\in \pa
K\subset\cE_{H}$. Finally, the local quasi-uniformity of the mesh
together with a scaling argument leads to the estimate
\eqref{local}.
\end{proof}

\begin{remark} An easy observation finds that the positive constant in \eqref{local} depends on the following ratio
\begin{equation}
\mu=\max\limits_{K\in\cT_H\backslash\cT_h}\max\limits_{\cT_h\ni T\subset K}\frac{h_K}{h_T}.
\end{equation}
In the analysis of  optimality of the adaptive finite element method,  this dependence is not allowed since we only know that $\cT_h$ is some refinement of $\cT_H$ by the newest vertex bisection
 and the boundness of $\mu$ is not guaranteed.
\end{remark}

To overcome the above difficulty, we  introduce a modified prolongation operator $J_h$ which preserves the local projection  property.   We need the prolongation operator $\Pi: W_H\rightarrow W_H^C$, where $W_H^C\subset W$ is
 some conforming finite element space over the mesh $\cT_H$.  Here we take
 $W_H^C$ as the Hsieh-Clough-Tocher finite element space over the mesh $\cT_H$ \cite{CiaBook,BrennerScott}.

 Let $\mathcal{F}$  be any (global) degree of freedom of $W_H^C$, i.e., $\mathcal{F}$ is either
the evaluation of a shape function or its first order derivatives at an interior node
of $\cT_H$, or the evaluation of the normal derivative of a shape function at a node on
an interior edge. For  $v_H\in W_H$, we  define \cite{BrennerSung2005}
\begin{equation}
\mathcal{F}(\Pi v_H)=\frac{1}{|\omega_{\mathcal{F}}|}\sum\limits_{K\in\omega_{\mathcal{F}}}\mathcal{F}(v_H|_K)
\end{equation}
where $\omega_{\mathcal{F}}$ is the set of triangles in $\cT_H$ that share the degree of freedom $\mathcal{F}$, and
 $|\omega_{\mathcal{F}}|$ is the number of elements of $\omega_{\mathcal{F}}$.  Then  a similar argument of \cite{BrennerSung2005} proves
 \begin{equation}\label{eq5.7}
 \|\na_H^2(v_H-\Pi v_H)\|_{L^2(\Om)}^2\lesssim \sum\limits_{e\in\cE_H}h_e\|[\na_H^2v_H\tau_e]\|_{L^2(e)}^2.
 \end{equation}
  Define
  \begin{equation*}
 \Omega_{\mathcal{R}}:=\textrm{interior}(\bigcup\{K: K\in\mathcal{T}_H\backslash\mathcal{T}_h,\}),
 \end{equation*}
 and
\begin{equation*}
 \Omega_{\mathcal{C}}:=\textrm{interior}(\bigcup\{K: K\in\mathcal{T}_H\cap\mathcal{T}_h, \partial K\cap \partial \Omega_{\mathcal{R}}=\emptyset \}).
 \end{equation*}
 The main idea herein is to take the mixture of the prolongation operators $I_{h}^{\prime}$
  and $\Pi$.  More precisely, we use $\Pi$ in the region $\Omega_{\mathcal{R}}$ where the elements of $\cT_H$ are refined and take $I_{h}^{\prime}$ on  $\Omega_{\mathcal{C}}$, and we define some mixture in the layer between them. This leads to the prolongation
   operator $J_{h}: W_H\rightarrow W_{h}$ by
   $$
   \mathcal{J}_{h}v_H=\left\{\begin{array}{ll}\Pi_{h}\Pi v_H &\text{ on }\Omega_{\mathcal{R}},\\ I_{h}^{\prime}v_H &\text{ on } \Omega_{\mathcal{C}},\\
 v_{h, tr}&\text{ on }\Omega\backslash(\Omega_{\mathcal{R}}\cup \Omega_{\mathcal{C}}),\end{array}\right.
   $$
 where $v_{h, tr}$ is defined by
   \begin{equation}
   \begin{split}
   v_{h, tr}(P)&=\left\{\begin{array}{ll}(\Pi v_H)(P) \text{ if  } P\in \pa\Om_{\mathcal{R}},\\(I_h^{\prime} v_H)(P) \text{ otherwise },\end{array}\right.\text{ for } P\in\mathcal{N}_h,\\
   \int_e\na_hv_{h, tr}\cdot\nu_e ds&=\left\{\begin{array}{ll}\int_e\na_h\Pi v_H\cdot\nu_eds & \text{ if } e\subset\pa\Omega_{\mathcal{R}}\\
      \int_e \na_hI_{h}^{\prime}v_H\cdot\nu_e ds &\text{ otherwise }\end{array}\right. \text{ for } e\in\cE_{h}.
   \end{split}
   \end{equation}

Define
$$
\mathcal{M}_{H,  h}:=\{K\in\cT_H, \partial K\cap \overline{\cup(\cT_H\backslash\cT_h)}\not=\emptyset \}.
$$

\begin{remark} It follows immediately from regularity of the mesh $\cT_H$ that
$$
\# \mathcal{M}_{H, h}\leq \kappa \# \cT_H\backslash\cT_h
$$
for a positive constant $\kappa\geq 1$ which is only dependent on the initial mesh $\cT_0$.
\end{remark}

\begin{lemma} It holds true that
\begin{equation}\label{eq3.5cb}
\|\nabla_{h}^2(J_h v_H-v_H)\|_{L^2(\Om)}^2 \lesssim
\sum\limits_{K\in\mathcal{M}_{H,  h}}\sum\limits_{e\subset\partial K}h_K\|[\nabla_H^2 v_H\tau_e]\|_{L^2(e)}^2\,~\text{for any}~ v_H\in W_H.
\end{equation}
\end{lemma}
\begin{proof} We only need to use the scaling argument like that in Lemma \ref{Lemma5.2} in the layer $\Omega\backslash(\Omega_{\mathcal{R}}\cup \Omega_{\mathcal{C}})$.  The desired result follows from the estimate \eqref{eq5.7} and the local projection property $J_hv_H|_K=v_H|_K$ for $K\in \Omega_{\mathcal{C}}$.
\end{proof}

\begin{lemma}\label{Lemma5.4}(Discrete reliability) It holds that
\begin{equation}\label{eq5.19}
\|u_{h}-u_H\|_{\cC_{h}}^2
 \lesssim \eta^2(u_H, \mathcal{M}_{H, h})\,.
\end{equation}
\end{lemma}
\begin{proof} For any $v_{h}\in W_{h}$,  we deduce from the discrete
problem \eqref{discrm} that
\begin{equation}\label{eq5.20}
\begin{split}
\|u_{h}-u_H\|_{\cC_{h}}^2=a_{h}(u_{h}-u_H,
u_{h}-v_{h})+a_{h}(u_{h}-u_H,v_{h}-u_H)
 =J_1+J_2,
\end{split}
\end{equation}
where
$$
J_1=\res_H(u_{h}-v_{h}), \text{ and }
J_2=a_{h}(u_{h}-u_H,v_{h}-u_H)\,.
$$
Thanks to  \eqref{eq2.10},  \eqref{eq3.4b} and \eqref{eq3.5b}, the
residual  $J_1$ can be bounded by a similar argument for the term on
the right hand-side of  \eqref{eq3.9}, which implies
\begin{equation}\label{eq5.21}
\begin{split}
 |J_1|&=|\res_H(u_{h}-v_{h})|=|\res_H((I-I_H)(u_{h}-v_{h}))|\\[0.5ex]
&\lesssim
\sum\limits_{K\in\mathcal{M}_{H, h}}h_K^2\|f\|_{L^2(K)}
\|\nabla_{h}^2(u_{h}-v_{h})\|_{L^2(K)}\,.
\end{split}
\end{equation}
 Since $v_h\in W_h$ is arbitrary, we apply the Young
and Cauchy-Schwarz inequalities in \eqref{eq5.20} to obtain that
\begin{equation}\label{eq5.22}
\begin{split}
 \|u_{h}-u_H\|_{\cC_{h}}^2
 \lesssim
 \sum\limits_{K\in\mathcal{M}_{H, h}}h_K^{4}\|f\|_{L^2(K)}^2
 +\inf\limits_{v_{h}\in
 W_{h}}\|v_{h}-u_H\|_{\cC_{h}}^2\,.
 \end{split}
\end{equation}
Taking $v_h=J_hu_H$ and applying  \eqref{eq3.5cb} 
complete the proof of the Lemma.
\end{proof}

We end this section by applying the previous discrete reliability to
show a result  indicating  that the bulk criterion is in some sense
a necessary condition for  reduction of the energy norm between
two levels.
\begin{lemma}\label{Lemma5.5}
If $\cT_{h}$ is a refinement of $\cT_H$ such that the following reduction holds
\begin{equation}\label{eq5.23}
\begin{split}
\|u-u_{h}\|_{\cC_{h}}^2+\osc^2(f,\cT_h)
 \leq
\al^{\prime} (\|u-u_{H}\|_{\cC_{H}}^2 +\osc^2(f,\cT_H))\,,
\end{split}
\end{equation}
for some $0<\al^{\prime}<1$, then there exists $0<\theta_{\ast}<1$
such that
\begin{equation}\label{eq5.24}
\theta_{\ast}\eta^2(u_H, \cT_H)\leq \eta^2(u_H,
\mathcal{M}_{H, h})\,.
\end{equation}
\end{lemma}
\begin{proof} We start with the following decomposition
\begin{equation}\label{eq5.25}
\begin{split}
&(1-\al^{\prime})(\|u-u_H\|_{\cC_{H}}^2+\osc^2(f,\cT_H))\\
 &\leq\|u-u_H\|_{\cC_{H}}^2+\osc^2(f,\cT_H)
-\|u-u_{h}\|_{\cC_{h}}^2-\osc^2(f,\cT_h)\\[0.5ex]
&=\|u_H-u_{h}\|_{\cC_{h}}^2
+2a_{h}(u-u_{h},u_h-u_H)+\osc^2(f,\cT_H)-\osc^2(f,\cT_h)\,.
\end{split}
\end{equation}
By the discrete reliability of Lemma \ref{Lemma5.4} with the
coefficient $C_{Drel}$, we have
\begin{equation}\label{eq5.26}
\|u_h-u_H\|_{\cC_{h}}^2 \leq C_{Drel}\eta^2(u_H,
\mathcal{M}_{H, h})\,.
\end{equation}
The quasi-orthogonality in Lemma \ref{Lemma3.2} with the coefficient
$C_{QO}$ yields
\begin{equation}
\begin{split}
 &|2a_{h}(u-u_{h},u_h-u_H)|\leq 2C_{QO}\|u-u_{h}\|_{\cC_h} \big(\sum\limits_{K\in
  \mathcal{M}_{H, h}}h_K^{4}\|f\|_{L^2(K)}^2\big)^{1/2}\,.
\end{split}
\end{equation}
It follows from \eqref{eq5.23} that
\begin{equation}
\|u-u_h\|_{\cC_h}\leq
\sqrt{\alpha^{\prime}}\big(\|u-u_H\|_{\cC_H}^2+\osc^2(f, \cT_H)
\big)^{1/2}.
\end{equation}
Therefore, we  apply the Cauchy-Schwarz inequality to obtain
\begin{equation}
\begin{split}
 |2a_{h}(u-u_{h},u_h-u_H)|&\leq
 \frac{1}{2}(1-\al^{\prime})\big(\|u-u_H\|_{\cC_H}^2+\osc^2(f,
 \cT_H)\big)\\[0.5ex]
&+2(C_{QO})^2\frac{\al^{\prime}}{1-\al^{\prime}}\sum\limits_{K\in
  \mathcal{M}_{H, h}}h_K^{4}\|f\|_{L^2(K)}^2\,.
\end{split}
\end{equation}
Since it is obvious that
 \begin{equation}\label{eq5.27}
 |\osc^2(f,\cT_H)-\osc^2(f, \cT_h)|\leq \eta^2(u_H,
 \mathcal{M}_{H, h})\,,
 \end{equation}
 we combine \eqref{eq5.25}- \eqref{eq5.27} to  prove the desired result by
 the parameter
 $$\theta_{\ast}=\frac{(1-\al^{\prime})^2C_{Eff}}{2(2\al^{\prime}(C_{QO})^2+(1-\al^{\prime})(C_{Drel}+1))}$$
 with  the efficiency constant $C_{Eff}$ of the estimator $\eta(u_H,
\cT_H)$ from Lemma \ref{Lemma2.1}.
\end{proof}

\section{Optimality}

To analyze the optimality, we follow an idea commonly used in the adaptive
finite element literature to introduce a nonlinear approximation class. First, we
have the following quasi-optimality.
\begin{equation}\label{eq5.28}
\|u-u_H\|_{\cC_H}^2
 \approxeq \inf\limits_{v_H\in W_H}\|u-v_H\|_{\cC_H}^2
  +\consis^2(u, \cT_H)\,,
\end{equation}
where the consistency error term is given by
\begin{equation}\label{eq5.29}
\consis(u, \cT_H)=\sup\limits_{v_H\in W_H}\frac{(f,
v_H)_{L^2(\Om)}-a_H(u, v_H)}
  {\|v_H\|_{\cC_H}}\,.
\end{equation}
It follows from \cite[Section 4.1]{Gudi2010} that
\begin{equation*}
\consis(u, \cT_H)\lesssim \inf\limits_{v_H\in W_H}\|u-v_H\|_{\cC_H}+\osc(f,\cT_H).
\end{equation*}
Therefore, we define
\begin{equation}\label{eq5.30}
\mathfrak{E}(N;u,
f):=\inf\limits_{\cT\in\mathbb{T}_N}\inf\limits_{v\in
W_{\cT}}\big(\|u-v\|_{\cC_{\cT}}^2+\osc^2(f,\cT)\big)\,.
\end{equation}
 Finally, we choose the nonlinear approximation class as follows:
\begin{equation}\label{eq5.31}
\mathbb{A}_s:=\big\{(u,f), |u, f|_s:=\sup\limits_{N>
N_0}N^{s}\mathfrak{E}(N;u, f)<+\infty\big\}.
\end{equation}
 Compared to  the adaptive conforming method for the second order
elliptic problem \cite{CasconKreuzerNochettoSiebert07,Stevenson06},
 we have not the
following monotone convergence:
\begin{equation}\label{monotone}
\inf\limits_{v_h\in W_h}\|u-v_h\|_{\cC_h}^2
  +\consis^2(u, \cT_h)^2\leq\inf\limits_{v_H\in W_H}\|u-v_H\|_{\cC_H}^2
  +\consis^2(u, \cT_H)\,,
\end{equation}
where $T_h$ is some refinement of $T_H$.  However, it follows from the quasi-orthogonality in Lemma \ref{Lemma3.2}, the efficiency
 of the estimator in Lemma \ref{Lemma2.1}, and the Young inequality that
\begin{equation}\label{eq5.32}
\|u-u_{h}\|_{\cC_{h}}^2+\osc^2(f,\cT_h)\leq C_2(\|u-u_{H}\|_{\cC_{H}}^2+\osc^2(f,\cT_H)).
\end{equation}

\begin{theorem}\label{Theorem5.6}
Let $\mathcal{M}_k$ be a set of marked elements with minimal
cardinality from Algorithm \ref{Algorithm},   $u$ the solution of
Problem \eqref{rmplate}, and  $(\cT_k, W_k, u_k)$  the sequence of
meshes, finite element spaces, and discrete solutions produced by
the adaptive finite-element methods with $0<\theta<
\frac{C_{Eff}}{2(2 (C_{QO})^2+C_{Drel}+1)}$. Then, the following
estimate holds:
\begin{equation}\label{eq5.33}
\#\mathcal{M}_k\lesssim (\al^{\prime})^{-\frac{1}{s}}|u,
f|_s^{\frac{1}{s}}(C_{2})^{\frac{1}{s}}\big(\|u-u_{k}\|_{\cC_k}^2+\osc^2(f,\cT_k)\big)^{-\frac{1}{s}}~
\text{for any}~ (u,f)\in\mathbb{A}_s,
\end{equation}
where the parameter $0<\al^{\prime}<1$ is from  Lemma
\ref{Lemma5.5}.
\end{theorem}
\begin{proof} We set $\epsilon=\al^{\prime} (C_{2})^{-1}
\big(\|u-u_k\|_{\cC_k}^2+\osc^2(f,\cT_k)\big)$  with
$0<\al^{\prime}<1$. Since $(u, f)\in\mathbb{A}_s$, there exists
 a $\cT_{\epsilon}$ of the refinement of $\cT_0$ and $u_{\epsilon} \in W_{\cT_{\epsilon}}$ with
 \begin{equation}\label{eq5.34}
 \begin{split}
 \#\cT_{\epsilon}-\#\cT_0\leq |u, f|_s^{1/s}\epsilon^{-1/s} \text{ and } \|u-u_{\epsilon}\|_{\cC_{\cT_{\epsilon}}}^2+\osc^2(f, \cT_{\epsilon})\leq \epsilon.
 \end{split}
 \end{equation}
   The overlay
 $\cT_{\ast}$ of $\cT_{\epsilon}$ and $\cT_k$ is the smallest refinement of
 both   $\cT_{\epsilon}$ and $\cT_k$.
 Let $u_{\ast}$ be the  finite element solution of \eqref{discrm} on the mesh
 $\cT_{\ast}$. Since $\cT_{\ast}$ is a refinement of
 $\cT_{\epsilon}$, we use,  \eqref{eq5.34},  and  \eqref{eq5.32} to
  obtain that
 \begin{equation}\label{eq5.35}
 \begin{split}
 \|u-u_{\ast}\|_{\cC_{\cT_{\ast}}}^2+\osc^2(f,\cT_{\ast})
 &\leq C_2(\|u-u_{\epsilon}\|_{\cC_{\cT_{\epsilon}}}^2
 +\osc^2(f,\cT_{\epsilon}))
 \\[0.5ex]
 &\leq
 C_{2}\epsilon=\al^{\prime}(\|u-u_{k}\|_{\cC_k}^2+\osc^2(f,\cT_k))\,.
 \end{split}
 \end{equation}
 We deduce from Lemma \ref{Lemma5.5}  that
 \begin{equation}\label{eq5.36}
\theta_{\ast} \eta^2(u_k, \cT_k)\leq
\eta^2(u_k, \mathcal{M}_{k, \ast}), \text{ for some }
\theta_{\ast}\in(0,1).
 \end{equation}
 We note that the step (3) in Algorithm \ref{Algorithm} with $\theta\leq \theta_{\ast}$  chooses  a subset of
  $M_k\subset\cT_k$ with minimal cardinality with the same property. Therefore
  \begin{equation}\label{eq5.37}
  \# \mathcal{M}_k\lesssim
  \#\cT_{\ast}-\#\cT_k\leq \#\cT_{\epsilon}-\#\cT_0.
  \end{equation}
 This together with the definition of $\epsilon$
leads to
  \begin{equation}\label{eq5.38}
\begin{split}
 \# \mathcal{M}_k
  \lesssim  (\al^{\prime})^{-\frac{1}{s}}|u, f|_{s}^{\frac{1}{s}}
  (C_{2})^{\frac{1}{s}}
  (\|u-u_{k}\|_{\cC_k}^2+\osc^2(f,\cT_k))^{-\frac{1}{s}},
  \end{split}
  \end{equation}
  which completes the proof.
 \end{proof}

\begin{theorem}\label{Theorem5.7}
Let the marking step in Algorithm \ref{Algorithm} select a set
$\mathcal{M}_k$ of marked elements with minimal cardinality,  $u$
the solution to Problem \eqref{discrm}, and  $(\cT_k, W_k, u_k)$ the
sequence of meshes, finite element spaces, and discrete solutions
produced by the adaptive finite-element methods with $0<\theta<
\frac{C_{Eff}}{2(2 (C_{QO})^2+C_{Drel}+1)}$. Then,
 it holds that
\begin{equation}\label{eq5.39}
\|u-u_N\|_{\cC_N}^2+\osc^2(f, \cT_N)\lesssim |u,
f|_s(\#\cT_N-\#\cT_0)^{-s}, \text{for}~(u, f)\in
\mathbb{A}_s.
\end{equation}
\end{theorem}
\begin{proof}  Let
$\mu=(\al^{\prime})^{-\frac{1}{s}}|u,
f|_{s}^{\frac{1}{s}}(C_{2})^{\frac{1}{s}}$.  We use  the result that
$\#\cT_k-\#\cT_0\lesssim \sum\limits_{j=0}^{k-1}\mathcal{M}_j$ from
\cite{Stevenson05,Stevenson06} to obtain that
\begin{equation}\label{eq5.40}
\begin{split}
&\#\cT_N-\#\cT_0\lesssim \sum\limits_{j=0}^{N-1}\mathcal{M}_j
\lesssim \mu\sum\limits_{j=0}^{N-1}(\|u-u_j\|_{\cC_j}^2+\osc^2(f,
\cT_j))^{-\frac{1}{s}}.
\end{split}
\end{equation}
It follows from the efficiency of the estimator that
\begin{equation}\label{eq5.42}
\|u-u_{j}\|_{\cC_j}^2+\osc^2(f,\cT_j)\approxeq \tilde{\eta}^2(u_j,
\cT_j),
\end{equation}
which gives
\begin{equation}\label{eq5.43}
\begin{split}
\|u-u_j\|_{\cC_j}^2+\gamma_1\tilde{\eta}^2(u_j,
 \cT_j) &\lesssim \|u-u_j\|_{\cC_j}^2+\osc^2(f,\cT_j).
\end{split}
\end{equation}
For any $ 0\leq j\leq N-1$, we use  the convergence result from Theorem
\ref{Theorem5.2} to   derive that
\begin{equation}\label{eq5.44}
\|u-u_N\|_{\cC_N}^2+ \gamma_1\tilde{\eta}^2(u_N, \cT_N)\leq
\al^{(N-j)}(\|u-u_j\|_{\cC_j}^2+\gamma_1\tilde{\eta}^2(u_j,
 \cT_j)).
\end{equation}
A combination of \eqref{eq5.40}-\eqref{eq5.44} yields
\begin{equation}\label{eq5.45}
\#\cT_N-\#\cT_0\lesssim \mu (\|u-u_N\|_{\cC_N}^2+\osc^2(f,
\cT_N))^{-1/s}\sum\limits_{j=1}^{N}\al^{j/s}.
\end{equation}
Setting $C_{\theta}=\al^{1/s}(1-\al^{1/s})^{-1},$   it is easy to prove that
\begin{equation}\label{eq5.46}
\sum\limits_{j=1}^{N}\al^{j/s}\leq C_{\theta}.
\end{equation}
Inserting this bound into \eqref{eq5.45} leads to
\begin{equation}\label{eq5.47}
\|u-u_N\|_{\cC_N}^2+\osc^2(f, \cT_N)\lesssim |u,
f|_s(\#\cT_N-\#\cT_0)^{-s},
\end{equation}
which completes the proof.
\end{proof}

\section{The extensions of the theory}
This section  extends  the theory to the Morley element in three
dimensions and the nonconforming linear elements in both two and
three dimensions.

\subsection{The Morley element in three dimensions}
Let $\cT_h$ be a decomposition of the domain $\Om\subset \R^3$ into
simplicies.  Given any face $F$, we let $\nu_F$ denote its unit
normal vector.  The Morley element in three dimensions  is defined
and analyzed in \cite{WX06}, where the space reads
\begin{equation}\label{threespace}
\begin{split}
&W_h:=\{v\in L^2(\Om), v|_K\in P_2(K), K\in\cT_h,  \int_e [v]\ds=0
\text{ for
any internal edge $e$},\\[0.5ex]
  &\int_e v\ds=0 \text{ for any boundary edge } e,
    \int_{F}[\na v\cdot \nu_F] dF=0 \text{ for any }\\[0.5ex]
    &\text{ internal face $F$}\,,
 \text{ and }\int_F \na v\cdot\nu_F dF=0 \text{ for any boundary
  face $F$}\}\,.
  \end{split}
\end{equation}

Define the estimator on each element $K\in\cT_h$ as
\begin{equation}
\eta_K=h_K^{2}\|f\|_{L^2(K)}+\bigg(\sum\limits_{F \subset\pa
K}h_F\|[\na_h^2u_h\times \nu_F]\|_{L^2(F)}^2 \bigg)^{1/2}\,,
\end{equation}
where $\times$ denotes the usual tensor product. The estimator  is defined by
\begin{equation}
\eta^2(u_h, \cT_h):=\sum\limits_{K\in \cT_h}\eta_K^2\,.
\end{equation}

The following reliability and efficiency of the estimator  were
proved in \cite{HuShi08}.

\begin{lemma} Let $u$ be the solution to the fourth  order elliptic problem with
 $u|_{\pa\Omega}=\frac{\pa u}{\pa\nu}|_{\pa\Omega}=0$
   in three dimensions,  $u_h$ be the finite element solution corresponding
   to the discrete space $W_h$ defined in \eqref{threespace}.  Then,
\begin{equation}
\|u-u_h\|_{\cC_h}\approxeq \eta_h
\end{equation}
up to the oscillation $\osc(f, \cT_h)$, where $\|\cdot\|_{\cC_h}$
and $\osc(f, \cT_h)$ are the three dimensional counterparts of the
discrete energy norm in \eqref{discreteenergy} and  the oscillation
in \eqref{eq2.8}, respectively.
 \end{lemma}

 \begin{lemma}\label{averagecon2}  Let $K_1, K_2\in\cT_h$ be two elements sharing a common face $F$
with three edges $e_{\ell}$ and midpoints $m_{\ell}$,
$\ell=1\,,2\,,3\,,$ and $v$ be a piecewise polynomial of degree
$\leq 1$ over $K_1\cup K_2$ such that
\begin{equation}\label{con2}
  v|_{K_1}(m_{\ell})=v|_{K_2}(m_{\ell}), \ell=1\,,2,3,  \text{ and } \int_{F}[\frac{\pa v}{\pa
\nu_{F}}]dF=0\,.
\end{equation}
 Then, $v$ is a polynomial of degree $\leq 1$  over $K_1\cup K_2$.
\end{lemma}

With these preparations, one can generalize the theories of the
quasi-orthogonality of Lemma \ref{Lemma3.2}, error reduction of
Theorem  \ref{Theorem5.2}, the discrete reliability of Lemma
\ref{Lemma5.4}, and the optimality of  Theorem  \ref{Theorem5.7} to
the Morley element method in three dimensions.

\subsection{The nonconforming linear elements for second order elliptic problems} In this subsection,
we let $\cT_h$ be a decomposition of the domain $\Om\subset \R^2$ or
$\Om\subset \R^3$ into simplicies in both two  and  three
dimensions. The nonconforming linear element spaces in both two and
three dimensions is defined by, respectively,
\begin{equation}\label{twospaceP1}
\begin{split}
&W_h:=\{v\in L^2(\Om), v|_K\in P_1(K), K\in\cT_h,  \int_e [v]\ds=0
\text{ for
any internal edge $e$},\\[0.5ex]
  &\int_e v\ds=0 \text{ for any boundary edge } e\}\,,
  \end{split}
\end{equation}
\begin{equation}\label{threespaceP1}
\begin{split}
&W_h:=\{v\in L^2(\Om), v|_K\in P_1(K), K\in\cT_h,  \int_F [v]\, dF=0
\text{ for
any internal face $F$},\\[0.5ex]
  &\int_F v\, dF=0 \text{ for any boundary face } F\}\,.
  \end{split}
\end{equation}
The continuous problems read: Given $f\in L^2(\Omega)$, find $u\in
H^1_0(\Omega)$ such that
\begin{equation}
(\nabla u,\nabla v)_{L^2(\Omega)}=(f, v)_{L^2(\Omega)} \text{ for
any }v\in H_0^1(\Omega).
\end{equation}
The discrete problems read: Given $f\in L^2(\Omega)$, find $u_h\in
W_h$ such that
\begin{equation}
(\nabla_h u_h,\nabla_h v_h)_{L^2(\Omega)}=(f, v_h)_{L^2(\Omega)}
\text{ for any }v_h\in W_h\,.
\end{equation}

The convergence of the adaptive nonconforming linear element methods
was first analyzed in \cite{CarstensenHoppe05b}. The theory in
Sections 3-7 can be extended to this case.  This
 extension gives another analysis of the convergence result from
 \cite{CarstensenHoppe05b}.

\section{Conclusion and comments}
In this paper, we carry out the convergence and optimality analysis
of the Morley element for the fourth order elliptic equation.
Moreover, we generalize the theory to the nonconforming linear
elements.  However, the analysis herein heavily depends on the
conservative properties of these two classes of nonconforming
elements  and the fact that the discrete stress is a piecewise
constant tensor. At the present time, it is unclear how  to generalize
these techniques to other nonconforming schemes of the fourth order
elliptic problems.

\end{document}